\documentclass[11pt]{amsart}

\usepackage{amsmath,amssymb,amsthm,amsfonts,mathtools}
\usepackage{geometry}
\usepackage{hyperref}
\usepackage{enumitem}
\usepackage{mathrsfs}
\usepackage{tikz-cd}
\usepackage{color}
\usepackage{verbatim}
\usepackage{bm}
\hypersetup{colorlinks=true, linkcolor=blue, citecolor=blue, urlcolor=blue}

\numberwithin{equation}{section}

\newtheorem{theorem}{Theorem}[section]
\newtheorem{proposition}[theorem]{Proposition}
\newtheorem{lemma}[theorem]{Lemma}
\newtheorem{corollary}[theorem]{Corollary}

\theoremstyle{definition}

\newtheorem{example}[theorem]{Example}

\newtheorem{remark}[theorem]{Remark}

\newcommand{\Hess}{\operatorname{Hess}}
\newcommand{\tr}{\operatorname{tr}}

\newcommand{\MCHess}{\operatorname{MCHess}}
\newcommand{\MHess}{\operatorname{MHess}}

\title[Hinge comparison and excess estimates via modified Hessian bounds]{Hinge comparison and excess estimates via modified Hessian bounds on manifolds with density}

\author{Nicholas Ng}

\address{1116 8th Avenue\\
Department of Mathematics and Statistics\\
Grinnell College\\
Grinnell, IA 50112}

\email{ngn11@grinnell.edu}
\urladdr{\url{https://www.nicholasngmath.com}}



\date{}

\begin{document}

\begin{abstract}
We study geometric consequences of modified Hessian estimates on Riemannian manifolds with density. For general radial functions, we derive tensorial, spectral, and trace identities for the modified Hessian in the original metric. We also establish a comparison identity relating the modified Hessians of two radial functions and characterize when they are proportional. We then prove a hinge comparison theorem under an upper bound on the modified Hessian, yielding a weighted law of cosines and lower bounds for the excess function. Finally, under nonnegative weighted sectional curvature and suitable density control assumptions, we obtain an Abresch–Gromoll-type upper excess estimate.
\end{abstract}

\maketitle

\section{Introduction}

Comparison geometry studies how curvature bounds control the behavior of
distance functions, geodesic spheres, Jacobi fields, and the geometry of
geodesic triangles. Classical sectional curvature bounds lead to Hessian
comparison, Laplacian comparison, and Toponogov's comparison theorem, which
provide quantitative information about the geometry of Riemannian manifolds
(see, e.g., \cite{CheegerEbin, Eschenburg, Lee18, Meyer, Petersen}).
Related comparison techniques under Ricci curvature bounds are discussed
in \cite{DaiWei}. Another important problem in comparison geometry is the
study of the excess function. Given three points $p,q,x\in M$, the excess is
defined by $$e(x)=d(p,x)+d(q,x)-d(p,q).$$ The classical Abresch--Gromoll excess estimate \cite{AbreschGromoll} gives an upper bound on the excess in terms of the distance from $x$ to a minimizing geodesic joining $p$ and $q$, together with the distance from $x$ to the endpoints. These estimates have played an important role in comparison geometry and in the study of the topology of manifolds with nonnegative Ricci curvature; see also \cite{DaiWei} for an exposition of comparison geometry and excess estimates.

Comparison geometry has also been developed in the setting
of Riemannian manifolds with density. The Bakry--Émery Ricci tensor was
introduced by Bakry--Émery \cite{BakryEmery}, and comparison geometry for
manifolds with density was developed by Wei--Wylie \cite{WeiWylie09}.
In particular, Wei--Wylie established comparison results under lower bounds
on the Bakry--Émery Ricci tensor, including weighted analogues of classical
comparison and excess estimates. A related but different notion of
curvature, called weighted sectional curvature, was introduced by Wylie
\cite{Wyl15} and further developed in subsequent work, including
\cite{KWY,WY16}. While the Bakry--Émery approach modifies the Ricci tensor to incorporate the density, weighted sectional curvature incorporates
the density into sectional curvature and is used to capture its
influence on the behavior of geodesics and Jacobi fields. This framework
uses a weighted connection and a corresponding conformal Hessian to
obtain comparison results for manifolds with density. In particular,
Kennard--Wylie--Yeroshkin \cite{KWY} developed comparison and rigidity
results in the setting of positive weighted sectional curvature.
For further details on the definition and geometric properties of
weighted sectional curvature, see Wylie \cite{Wyl15} and
Wylie--Yeroshkin \cite{WY16}. For additional background and developments
in the geometry of manifolds with density, see \cite{KWY, Wyl15,WY16}.

The present work continues our study of comparison geometry on Riemannian manifolds $(M, g)$ with density $\varphi$ through the modified Hessian tensor \begin{equation}\label{eq:1.1}\MHess u=\Hess u-\langle\nabla u,\nabla\varphi\rangle g,
\end{equation} where $u \in C^\infty(M)$ is a smooth function. In particular, we study radial functions of the form $u = h \circ r_p$, where $r_p = d(p, \cdot)$ is the distance function from a fixed point $p \in M$ and $h$ is a smooth function of a single variable. The tensor \eqref{eq:1.1} is the expression of the modified
conformal Hessian introduced by Kennard--Wylie--Yeroshkin \cite{KWY} in the original metric $g$.
In \cite{NgPaper1}, we obtained upper bounds on the modified Hessian of
the square of a radial distance function under nonnegative weighted
sectional curvature and suitable density control assumptions. These bounds
yield Hessian comparison, shape operator comparison, and weighted Laplacian
comparison, together with consequences for radial volume density and
volume growth. The purpose of the present paper is to obtain further
geometric consequences of these modified Hessian estimates. In particular,
we study comparison estimates for geodesic hinges and derive lower and
upper bounds on the excess function.

Our first result concerns the structural properties of the modified
Hessian for general radial functions. We establish a radial-tangential
decomposition in terms of the modified shape operator
$S_\varphi=S-\varphi_r I$, with the radial component determined by
$h''-h'\varphi_r$ and the tangential component given by $h'S_\varphi$.
We also obtain the corresponding trace identity, relating the trace
of the modified Hessian to the weighted Laplacian and the modified
shape operator. These identities hold on the regular set of the
distance function and require no curvature assumptions. While the
radial-tangential decomposition is implicit in the modified conformal
Hessian framework of Kennard--Wylie--Yeroshkin \cite{KWY}, our
formulation works directly with the original metric $g$ and makes
explicit the resulting tensorial and trace identities.

We then obtain a hinge comparison theorem for radial functions. Let $p,q,x\in M$, and let $\sigma:[0,c]\to M$ be a unit-speed minimizing
geodesic from $q$ to $x$. Set $$a=d(p,q),\qquad b=d(p,x),\qquad c=d(q,x),$$
and let $\beta$ denote the interior angle at $q$ determined by the directions $-\nabla r_p(q)$ and $\sigma'(0)$, where $r_p=d(p,\cdot)$. Suppose that $u=h\circ r_p$ and that $$\MHess u\leq K(r_p)g$$ along $\sigma$, for some continuous function $K$. We obtain the following inequality: $$h(b)\leq h(a)-ch'(a)\cos\beta
+\int_0^c(c-s)
\left[K(R(s))+h'(R(s))\varphi_r(\sigma(s))\right]\,ds,$$
where $R(s)=r_p(\sigma(s))$ and $$\varphi_r=\langle\nabla\varphi,\nabla r_p\rangle.$$ This gives an integral comparison along the side of a geodesic hinge. Unlike the classical hinge comparison in Toponogov's theorem \cite{CheegerEbin, Meyer}, the comparison does not involve a corresponding hinge in a model space. Instead, it follows directly from an upper bound on the modified Hessian of a radial function based at $p$.

Applying the hinge comparison theorem to the quadratic radial function $u=\frac{1}{2}r_p^2$, we obtain a weighted law of cosines. In particular, under a uniform upper bound on $\MHess(\frac{1}{2}r_p^2)$ and suitable radial density control, we obtain $$b^2\leq a^2-2ac\cos\beta+(K_0+c_0)c^2.$$ We then use this inequality to derive lower bounds on the excess $$e=a+c-b.$$ These estimates give quantitative control of the excess in terms of the side lengths and the interior angle of a geodesic hinge. In the Euclidean case with constant density, the weighted law of cosines reduces to the classical law of cosines, and the resulting excess inequality recovers the corresponding Euclidean excess identity.

We also obtain an upper excess estimate of Abresch--Gromoll type. Let
$(M^n,g,\varphi)$ be a complete Riemannian manifold with density, with
$n\geq3$, and suppose that
$$\overline{\sec}_{\varphi}\geq0,\qquad
a_0\leq\varphi\leq0,\qquad
|\nabla\varphi|\leq\frac{c_0}{r_p},$$ where $r_p=d(p,\cdot)$. Given a minimizing geodesic $\gamma$ from $p$ to $q$, let $$h=d(x,\gamma),\qquad
s=\min\{d(p,x),d(q,x)\}.$$ Here $h$ measures how far $x$ lies from $\gamma$ and is often referred to as the height of $x$ above $\gamma$. We prove that, whenever $h\leq s/2$, $$e(x)\leq\mathcal C_A
\left(\frac{h^{A+1}}{s}\right)^{1/A},$$ where $$A=(n-1)e^{-2a_0}+2c_0 \qquad \text{and} \qquad \mathcal C_A = \frac{2A}{A-1}\left(\frac{2A}{A+1}\right)^{1/A}.$$ In particular, the exponent in this estimate depends on the dimension as well as the bounds on the density and its gradient. In the case $\varphi = 0$, $a_0=c_0=0$, and hence $A=n-1$. Thus, we recover the exponent $1/(n-1)$ in the classical Abresch--Gromoll estimate \cite{AbreschGromoll}. We note, however, that the curvature assumption in the zero-density case
is simply nonnegative sectional curvature, which is stronger than the nonnegative Ricci curvature assumption in the classical result. The estimate therefore provides an excess comparison under nonnegative weighted sectional curvature, using the modified Hessian and weighted Laplacian comparison developed in \cite{NgPaper1}.

Hurtado--Palmer--Rosales \cite{HurtadoPalmerRosales} establish Hessian comparison results under radial lower bounds on certain weighted sectional curvatures, which are related to those introduced by Wylie. Our approach instead uses the modified Hessian, $\MHess u$, associated with Wylie's weighted sectional curvature framework to derive structural identities, hinge comparison, and
excess estimates. Park \cite{Park2025} studies Hessian estimates for the Green function and derives associated
monotonicity formulas, whereas our work concerns the modified Hessian for general radial functions, including squared distance functions. We also note that Wei--Wylie \cite{WeiWylie09} obtain a weighted version of the Abresch--Gromoll inequality under lower Bakry--Émery Ricci curvature bounds and suitable bounds on the density, while Gigli--Mosconi \cite{GigliMosconi} establish the Abresch--Gromoll inequality in the non-smooth setting of $\mathrm{CD}(K,N)$ spaces.

The paper is organized as follows. In Section \ref{sec:2}, we establish structural identities for the modified Hessian of radial functions. In particular, we obtain a radial-tangential decomposition and a comparison identity
relating the modified Hessians of two radial functions. In Section \ref{sec:3}, we prove the hinge comparison theorem and derive the weighted law of cosines. Section \ref{sec:4} is devoted to lower excess estimates and their consequences for the geometry of geodesic triangles. Finally, in Section \ref{sec:5}, we establish the Abresch--Gromoll-type upper excess estimate under nonnegative weighted sectional curvature and density control assumptions.

\section{Preliminaries}\label{sec:2}

\subsection{Weighted sectional curvature}

In \cite{Wyl15}, Wylie introduced the notion of \textit{weighted sectional curvature} for a Riemannian manifold $(M,g)$ with density $\varphi$, given by $$\overline{\sec}_\varphi(U,V) = \sec(U,V)+\Hess\varphi(U,U)+d\varphi(U)^2 = g(R^{\nabla^\varphi}(V,U)U,V),$$
for orthonormal vectors $U,V$, where $\nabla^\varphi$ is the \textit{weighted connection}
$$\nabla_X^\varphi Y = \nabla_XY-d\varphi(X)Y-d\varphi(Y)X,$$
and $R^{\nabla^\varphi}$ denotes the corresponding weighted curvature tensor
$$R^{\nabla^\varphi}(X,Y)Z = \nabla_X^\varphi\nabla_Y^\varphi Z
-\nabla_Y^\varphi\nabla_X^\varphi Z
-\nabla_{[X,Y]}^\varphi Z.$$
Here $\nabla$ denotes the Levi-Civita connection of the metric $g$. Comparison geometry for weighted sectional curvature was further developed by Kennard--Wylie--Yeroshkin in \cite{KWY}. For further details on the definition and geometric properties of weighted sectional curvature, see \cite{KWY, Wyl15,WY16}.

Throughout this paper, we use the weighted measure
$$d\mu=e^{-(n+1)\varphi}\,d\operatorname{vol}_g,$$ The factor $n+1$ is chosen to agree with the normalization of the weighted connection in \cite{KWY}.
The corresponding weighted Laplacian is defined by
$$\Delta_\mu f = \Delta f-(n+1)\langle\nabla\varphi,\nabla f\rangle.$$

\subsection{Modified Hessians}

Given a manifold with density $(M,g,\varphi)$, consider the conformal metric
$$\widetilde g=e^{-2\varphi}g.$$
It is well known that for any smooth function $u$,
\begin{equation}\label{eq:2.1}
\Hess_{\widetilde g}u = \Hess_g u + d\varphi\otimes du + du\otimes d\varphi -\langle\nabla\varphi,\nabla u\rangle g.\end{equation}
A basic property in the unweighted setting is that away from the cut locus, the radial direction $\nabla r$ lies in the kernel of $\Hess_g r$ whenever $r$ is a distance function. However, this is not generally true for $\Hess_{\widetilde g}r$. To remedy this, Kennard--Wylie--Yeroshkin in \cite{KWY} considered the following lower order perturbation of $\Hess_{\widetilde g}u$:
$$
\Hess_{\widetilde g}u
-d\varphi\otimes du
-du\otimes d\varphi.
$$
We set notation for the \textit{modified conformal Hessian}:
$$
\MCHess u
:=
\Hess_{\widetilde g}u
-d\varphi\otimes du
-du\otimes d\varphi,
$$
where $\widetilde g=e^{-2\varphi}g$ and $u$ is any smooth function. Substituting \eqref{eq:2.1} in the conformal Hessian term of $\MCHess u$ gives the \textit{modified Hessian}
$$
\MHess u
:=
\Hess_g u-\langle\nabla u,\nabla\varphi\rangle g.
$$
The modified Hessian framework considered here is motivated by the modified conformal Hessian tensor introduced by Kennard--Wylie--Yeroshkin in \cite{KWY}.

\begin{remark}
All geometric quantities in this paper are computed with respect to the underlying metric $g$. In particular, geodesics, distance functions, Hessians, and curvature quantities are those of the metric $g$.
\end{remark}

\subsection{Hessian identities and radial distance functions}\label{subsec:2.3}

Fix a point $p\in M$ and let
$$
r=r_p=d(p,\cdot).
$$
We denote the regular set of the distance function by
$$
\mathcal R_p=M\setminus\bigl(\{p\}\cup\operatorname{Cut}(p)\bigr).
$$
Wherever the distance function $r$ is smooth, we have
$$
|\nabla r|=1
\qquad\text{and}\qquad
\Hess r(\nabla r,\cdot)=0.
$$
We also define
$$
\varphi_r:=\langle\nabla\varphi,\nabla r\rangle,
\qquad
g_r:=g-dr\otimes dr.
$$
Here $\varphi_r$ denotes the radial derivative of the density function, and $g_r$ is the restriction of the metric to directions perpendicular to $\nabla r$, extended by zero in the radial direction. In geodesic polar coordinates, if $\varphi=\varphi(r,\theta)$, then $\varphi_r=\partial_r\varphi(r,\theta)$. Unless otherwise stated, we do not assume that the density function $\varphi$ is globally radial.

More generally, let $h$ be a smooth function of one variable and consider the radial function
$$u=h\circ r.$$
Wherever $r$ is smooth, the gradient and Hessian of $u$ are given by
$$\nabla u=h'(r)\nabla r$$
and
$$\Hess u = h'(r)\Hess r+h''(r)dr\otimes dr.$$
In particular, for the quadratic radial function
$$u(r)=\frac12r^2,$$
we have
$$\nabla u=r\nabla r, \qquad \Hess u=dr\otimes dr+r\Hess r.$$

\subsection{Shape operator}

Let $$S_r(p)=\{x\in M:d(p,x)=r\}$$
denote the geodesic sphere of radius $r$ centered at $p$. Wherever the distance function $r$ is smooth, define the shape operator by
$$S(v)=\nabla_v\nabla r, \qquad v\in(\nabla r)^\perp.$$
Then
$$\Hess r(v,w)=g(Sv,w)$$
for $v,w\in(\nabla r)^\perp$. In particular, $S$ is self-adjoint with respect to the induced metric on the tangent spaces of the geodesic spheres. See \cite{Lee18,Petersen} for background on the shape operator and the geometry of geodesic spheres.

\subsection{Modified Hessian Tensor Identities}

The modified Hessian tensor admits a natural decomposition into radial
and tangential components. Kennard--Wylie--Yeroshkin derived this radial-tangential decomposition in their study of weighted convexity \cite{KWY}, using the modified conformal Hessian to obtain geometric consequences of weighted sectional curvature bounds. Here, we build on their framework by investigating the structural properties of the modified Hessian in the original metric $g$, making explicit its tensorial, spectral, and trace properties rather than working within the conformal metric $\tilde g$ and setting of weighted convexity.

For $X, Y \in (\nabla r)^\perp$, define $S_\varphi$ by $$g(S_\varphi X, Y) = (\Hess r - \varphi_r g_r)(X, Y).$$ Then $S_\varphi = S - \varphi_r I$, where $S$ is the usual shape operator defined by $g(SX,Y) = \Hess r(X,Y)$ on $(\nabla r)^\perp$. We call $S_\varphi$ the \textit{modified shape operator}.

\begin{proposition}[Properties of the modified Hessian]
   Let $D \subseteq \mathcal R_p$ be an open subset, let $I \subset (0,\infty)$ be an interval containing $r(D)$, and let $h:I \rightarrow \mathbb R$ be a smooth radial function and let $u = h \circ r$. Then the following holds.

    \begin{enumerate}
        \item (Radial and Tangential decomposition) The modified Hessian can be decomposed into radial and tangential components, $$\MHess u = h'(\Hess r - \varphi_rg_r) + (h'' - h'\varphi_r)dr \otimes dr,$$ with respect to the orthogonal decomposition $TM\big|_D = (\nabla r)^\perp \oplus \mathrm{Span}\{\nabla r\}$.

        \item (Radial direction) $$\MHess u (\nabla r, \nabla r) = h'' - h'\varphi_r.$$

        \item (Tangential direction) For every $X, Y \in (\nabla r)^\perp$, $$\MHess u (X, Y) = h'(\Hess r - \varphi_rg_r)(X, Y).$$

        \item (Orthogonal directions) For any $X \perp \nabla r$, $$\MHess u (\nabla r, X) = 0.$$

        \item (Trace and weighted Laplacian identities) The weighted Laplacian of a radial potential $u$ with respect to the measure $\mu = e^{-(n + 1)\varphi}d\mathrm{vol}_g$ is $$\Delta_\mu u = h'' + h'\Delta_\mu r.$$ Thus \begin{align*}
\tr_g(\MHess u)
&=\Delta u-n\left<\nabla\varphi,\nabla u\right>\\
&=\Delta_\mu u+h'\varphi_r\\
&=h''+h'(\Delta r-n\varphi_r)\\
&=h'\tr(S_\varphi)+h''-h'\varphi_r.
\end{align*}
    \end{enumerate}
\end{proposition}

\begin{proof}
Since $u=h\circ r$, we have
$$\nabla u=h'\nabla r$$
and
$$\Hess u=h'\Hess r+h''dr\otimes dr.$$
By the definition of the modified Hessian,
$$\MHess u=\Hess u-\left<\nabla u,\nabla\varphi\right>g,$$
and since $\left<\nabla u,\nabla\varphi\right>=h'\varphi_r$, it follows that
\begin{align*}
\MHess u
&=h'\Hess r+h''dr\otimes dr-h'\varphi_r g\\
&=h'(\Hess r-\varphi_r g_r)
+(h''-h'\varphi_r)dr\otimes dr,
\end{align*}
where we have used $g=g_r+dr\otimes dr$. This proves (1).

Since $|\nabla r|=1$ on $\mathcal R_p$, we have
$$dr(\nabla r)=1$$
and
$$\Hess r(\nabla r,\cdot)=0.$$
Moreover, $g_r(\nabla r,\cdot)=0$. Evaluating the decomposition in the radial direction gives
$$\MHess u(\nabla r,\nabla r)=h''-h'\varphi_r,$$
which proves (2).

For $X,Y\in(\nabla r)^\perp$, we have $dr(X)=dr(Y)=0$. Thus
$$\MHess u(X,Y)=h'(\Hess r-\varphi_r g_r)(X,Y),$$
proving (3).

Similarly, for $X\perp\nabla r$, we have $dr(X)=0$ and
$$(\Hess r-\varphi_r g_r)(\nabla r,X)=0.$$
Consequently,
$$\MHess u(\nabla r,X)=0,$$
which proves (4).

Finally, for (5), the weighted Laplacian associated with the measure
$d\mu=e^{-(n+1)\varphi}d\mathrm{vol}_g$ is
$$\Delta_\mu f=\Delta f-(n+1)\left<\nabla\varphi,\nabla f\right>.$$
Using $\nabla u=h'\nabla r$ and $\Delta u=h'\Delta r+h''$, we obtain
\begin{align*}
\Delta_\mu u
&=\Delta u-(n+1)\left<\nabla\varphi,\nabla u\right>\\
&=h'\Delta r+h''-(n+1)h'\varphi_r\\
&=h''+h'\Delta_\mu r.
\end{align*}
Moreover, tracing the definition of the modified Hessian gives
\begin{align*}\tr_g(\MHess u ) &= \Delta u - n\left<\nabla \varphi, \nabla u\right> \\ &= \Delta u - (n + 1)h'\varphi_r + h'\varphi_r \\ &= \Delta_\mu u +h'\varphi_r \\ &= h'' + h'\Delta_\mu r + h'\varphi_r \\ &= h'' + h'(\Delta r - (n + 1)\varphi_r) + h'\varphi_r \\ &= h'(\Delta r - (n-1)\varphi_r) - h'\varphi_r + h'' \\ &=h'\tr(S_\varphi) + h'' - h'\varphi_r.\end{align*}
\end{proof}

\begin{remark}
The modified Hessian is a density-dependent modification of the Hessian tensor whose trace is related, but not identical, to the weighted Laplacian. Since the term ``weighted Hessian'' might suggest that its trace coincides with the weighted Laplacian, we use the term \textit{modified Hessian} instead. The term \textit{density-modified
Hessian} would also be appropriate, but we adopt the shorter terminology.
\end{remark}

In addition to the $(0, 2)$-tensor $\MHess u$, we also use $\MHess u$ to denote the corresponding $(1, 1)$-tensor (or endomorphism). Since the $(0, 2)$-tensor is symmetric, the corresponding endomorphism is self-adjoint.

\begin{proposition}
 Let $D \subseteq \mathcal R_p$ be an open subset, let $h:I \rightarrow \mathbb R$ be a smooth radial function on an interval containing $r(D)$, and let $u = h \circ r$. \begin{enumerate}
        \item The orthogonal decomposition of the tangent bundle $$TM\big|_D = \mathrm{Span}\{\nabla r \} \oplus (\nabla r)^\perp$$ is invariant under $\MHess u$.

        \item With respect to any orthonormal frame $\{E_1, ... , E_{n - 1}, \nabla r\}$, $\MHess u$ has block diagonal representation $$\MHess u = \begin{bmatrix}
            h'[S_\varphi] & 0 \\ 0 & h'' - h'\varphi_r
        \end{bmatrix},$$ where $[S_\varphi]$ is the matrix of the modified shape operator with respect to the orthonormal basis $\{E_1, ..., E_{n -1 }\}$. 

        \item Since $S_\varphi$ is self-adjoint, at each point there exists an orthonormal basis $\{E_1, ... , E_{n -1}, \nabla r\}$ such that $S_\varphi E_i = \lambda_i E_i$ for $i = 1, ..., n -1$; we refer to the eigenvalues $\lambda_i$ of $S_\varphi$ as ``modified principal curvatures'' and their corresponding eigenvectors $E_i$ as ``modified principal directions''. Consequently, the modified Hessian of $u$ is a diagonal matrix $$\MHess u = \mathrm{Diag}\left(h'\lambda_1, ... , h'\lambda_{n -1 }, h'' - h' \varphi_r\right),$$ and the spectrum of $\MHess u$ is $$\mathrm{Spec}(\MHess u ) = \{h'\lambda_1, ..., h'\lambda_{n - 1}, h'' - h'\varphi_r\}.$$ For the quadratic modified distance function, $$\mathrm{Spec}\left(\MHess\left(\dfrac{1}{2}r^2\right)\right) = \{r\lambda_1, ..., r\lambda_{n - 1}, 1 - r\varphi_r\}.$$
    \end{enumerate}
\end{proposition}

\begin{proof}
By the radial-tangential decomposition established in the preceding proposition, we have
$$\MHess u(\nabla r,X)=0$$
for every $X\perp\nabla r$. Since $\MHess u$ is self-adjoint, it follows that both $\mathrm{Span}\{\nabla r\}$ and $(\nabla r)^\perp$ are invariant under the corresponding endomorphism. This proves (1).

On the tangential subspace $(\nabla r)^\perp$, the modified Hessian acts as $h'S_\varphi$, while on the radial direction it acts by multiplication by $h''-h'\varphi_r$. Consequently, with respect to the stated orthonormal frame, we obtain
$$
\MHess u=
\begin{bmatrix}
h'[S_\varphi] & 0\\
0 & h''-h'\varphi_r
\end{bmatrix},
$$
which proves (2).

Finally, since $S_\varphi$ is self-adjoint on $(\nabla r)^\perp$, the spectral theorem guarantees the existence of an orthonormal basis $\{E_1,\ldots,E_{n-1}\}$ such that
$$S_\varphi E_i=\lambda_i E_i$$
for $i=1,\ldots,n-1$. Combining this basis with $\nabla r$ gives
$$
\MHess u=\mathrm{Diag}(h'\lambda_1,\ldots,h'\lambda_{n-1},h''-h'\varphi_r).
$$
Therefore,
$$
\mathrm{Spec}(\MHess u)
=\{h'\lambda_1,\ldots,h'\lambda_{n-1},h''-h'\varphi_r\}.
$$
Taking $h(r)=\frac12r^2$, so that $h'=r$ and $h''=1$, gives
$$
\mathrm{Spec}\left(\MHess\left(\frac12r^2\right)\right)
=\{r\lambda_1,\ldots,r\lambda_{n-1},1-r\varphi_r\},
$$
as desired.
\end{proof}

\begin{theorem}[Comparison identity for modified Hessians] Let $D \subseteq \mathcal R_p$ be an open subset, and let $I \subseteq (0, \infty)$ be a connected interval containing $r(D)$. Given two smooth radial functions $u = h \circ r$ and $u_0 = h_0 \circ r$ where $h, h_0 \in C^2(I)$ and $h_0' \ne 0$ on $I$, on $D$ we have $$\MHess (h \circ r) = \dfrac{h'(r)}{h_0'(r)}\MHess (h_0 \circ r) + \left(h''(r) - \frac{h'(r)}{h_0'(r)}h_0''(r)\right)dr \otimes dr.$$
\end{theorem}

\begin{proof}
By the radial-tangential decomposition established above, we have
$$\MHess(h\circ r) = h'(\Hess r-\varphi_r g_r) +(h''-h'\varphi_r)dr\otimes dr$$
and $$\MHess(h_0\circ r) = h_0'(\Hess r-\varphi_r g_r) +(h_0''-h_0'\varphi_r)dr\otimes dr.$$
Multiplying the second identity by $\dfrac{h'}{h_0'}$, we obtain $$\frac{h'}{h_0'}\MHess(h_0\circ r) = h'(\Hess r-\varphi_r g_r) + \left(\frac{h'}{h_0'}h_0''-h'\varphi_r\right)
dr\otimes dr.$$ Subtracting this expression from the first identity gives
\begin{align*}
\MHess(h\circ r)
-\frac{h'}{h_0'}\MHess(h_0\circ r)
&=\left(h''-\frac{h'}{h_0'}h_0''\right)
dr\otimes dr.
\end{align*}
Rearranging establishes the desired identity.
\end{proof}

The preceding identity shows that the modified Hessian of any radial potential can be expressed in terms of the modified Hessian of another radial potential, up to a scalar rescaling and the addition of a radial term. It is natural to ask when the radial term vanishes. Note that $$h'' - \frac{h'}{h_0'}h_0'' = 0 \iff \left(\frac{h'}{h_0'}\right)' = 0.$$ Integrating twice gives $$h = Ch_0 + C_1,$$ for some integration constants $C$ and $C_1$. In this case, one radial function is affine in the other and $$\MHess u = C\MHess u_0.$$ This motivates the following proposition.

\begin{proposition}[Characterization of proportional modified Hessians] Let $D \subseteq \mathcal R_p$ be a connected open subset, $I \subseteq (0, \infty)$ a connected interval containing $r(D)$, and $h, h_0 \in C^2(I)$ with $h_0' \ne 0$ on $I$. Let $u = h \circ r$ and $u_0 = h_0 \circ r$. \begin{enumerate}
    \item If $h = Ch_0 + C_1$ on $I$, then $\MHess u = C \MHess u_0$. 
    
    \item Suppose $\MHess u = \Lambda \MHess u_0$ on $D$ where $\Lambda \in C^{\infty}(D)$. If $\Hess r - \varphi_rg_r \ne 0$ for all $p \in D$, then $h = Ch_0 + C_1$ on $r(D)$ for some constants $C, C_1$, where $\Lambda \equiv C$. 

    \item If $\Hess r = \varphi_r g_r$ on $D$, then $$\MHess u = \Lambda \MHess u_0$$ is equivalent to $$h'' - h'\varphi_r = \Lambda(h_0'' - h_0'\varphi_r),$$ where the proportionality factor $\Lambda$ need not be constant.
\end{enumerate}

\begin{proof}
For (1), suppose that $h=Ch_0+C_1$. Then $h'=Ch_0'$ and $h''=Ch_0''$. Substituting into the comparison identity, we obtain
\begin{align*}
\MHess u
&=\frac{h'}{h_0'}\MHess u_0
+\left(h''-\frac{h'}{h_0'}h_0''\right)dr\otimes dr\\
&=C\MHess u_0+(Ch_0''-Ch_0'')dr\otimes dr\\
&=C\MHess u_0,
\end{align*}
as desired.

For (2), suppose that $\MHess u=\Lambda\MHess u_0$ on $D$. By the comparison identity,
\begin{align*}
\left(\frac{h'}{h_0'}-\Lambda\right)\MHess u_0
+\left(h''-\frac{h'}{h_0'}h_0''\right)
dr\otimes dr=0.
\end{align*}
Restricting this identity to tangential vectors $X,Y\in(\nabla r)^\perp$, we obtain
$$
\left(\frac{h'}{h_0'}-\Lambda\right)
h_0'(\Hess r-\varphi_r g_r)(X,Y)=0.
$$
Since $h_0'\ne0$ and $\Hess r-\varphi_r g_r$ is nonzero at every point of $D$, it follows that
$$\Lambda=\frac{h'}{h_0'}$$
throughout $D$. Substituting this into the radial component of the comparison identity gives
$$h''-\frac{h'}{h_0'}h_0''=0.$$
Equivalently,
$$\left(\frac{h'}{h_0'}\right)'=0.$$
Since $D$ is connected and $r$ is smooth, $r(D)$ is a connected interval. Therefore, there exists a constant $C$ such that
$$\frac{h'}{h_0'}=C$$
on $r(D)$. Integrating, we obtain
$$h=Ch_0+C_1$$
on $r(D)$ for some constant $C_1$. Consequently, $\Lambda\equiv C$ on $D$.

For (3), suppose that $\Hess r=\varphi_r g_r$ on $D$. Then the radial-tangential decomposition gives
$$\MHess u=(h''-h'\varphi_r)dr\otimes dr$$
and
$$\MHess u_0=(h_0''-h_0'\varphi_r)dr\otimes dr.$$
Therefore,
$$\MHess u=\Lambda\MHess u_0$$
if and only if
$$h''-h'\varphi_r
=\Lambda(h_0''-h_0'\varphi_r),$$
as desired.
\end{proof}
    
\end{proposition}

\begin{remark} Wherever $r_p$ is smooth, the modified Hessian of any radial function $u = h\circ r$ can be written in terms of the quadratic modified distance function $u_0(r) = \frac{1}{2}r^2$: $$\MHess u = \frac{h'}{r}\MHess \left(\frac{1}{2}r^2\right) + \left(h'' - \frac{h'}{r}\right)dr \otimes dr.$$ Moreover, if the modified Hessian of the quadratic distance function is uniformly bounded above on a region $D$ where $r$ is smooth, i.e., if $\MHess(\frac{1}{2}r^2) \leq Kg$ for some constant $K \geq 0$, and \begin{center}$0 \leq \dfrac{h'}{r} \leq P$ \, \, and \, \,  $h'' - \dfrac{h'}{r} \leq Q$ \, \, \text{on $D$},\end{center} then $$\MHess u \leq (PK + \max\{Q, 0\})g.$$ Note that the sign of $h'/r$ is needed for such a bound, but there is no sign restriction on $Q$. More generally, $$\|\MHess u \| \leq \frac{|h'|}{r}\left\|\MHess\left(\frac{1}{2}r^2\right)\right\| + \left|h'' - \frac{h'}{r}\right|$$ where $\|A\| = \|A\|_{\mathrm{op}} = \sup_{|X| = 1}|AX|$ so that $\|A\| = \max_i |\lambda_i|$ for self-adjoint operators $A$ and their eigenvalues $\lambda_i$. Thus if $\|\MHess((1/2)r^2)\| \leq K_0$ on a domain $D$ where $r$ is smooth, and \begin{center}
$\sup\limits_D \dfrac{|h'|}{r} \leq P_0$ \, \, and \, \, $\sup\limits_D \left|h'' - \dfrac{h'}{r}\right| \leq Q_0$ \, \,  where $P_0, Q_0 \geq 0$,
\end{center} then $$\|\MHess u\| \leq P_0K_0 + Q_0.$$\end{remark}

\section{Weighted Hinge Comparison}\label{sec:3}

We use the modified Hessian of a radial function to derive an integral comparison along a minimizing geodesic forming one side of a hinge.

\begin{theorem}[Weighted Hinge Comparison]\label{thm:3.1}

Let $(M^n, g, \varphi)$ be an $n$-dimensional Riemannian manifold with density, and fix a basepoint $p \in M$. Let $r_p(x) = d(p, x)$, and suppose $u = h \circ r_p$ where $h \in C^2((0, \infty))$. Let $K:(0, \infty) \rightarrow \mathbb R$ be continuous and let $\sigma:[0, c] \rightarrow M$ be a unit-speed minimizing geodesic from $q$ to $x$ where $r_p$ is smooth on a neighborhood of $\sigma([0, c])$. Set \begin{center}
    $a = r_p(q) = d(p, q), \, \, \, \, \, b = r_p(x) = d(p, x), \, \, \, \, \, \text{and} \, \, \, \, \, R(s) = r_p(\sigma(s)),$
\end{center} where $p \ne q$ and $q \ne x$. Suppose $\MHess u \leq K(r_p)g$ along $\sigma$. Then $$u(x) \leq u(q) + c(u\circ\sigma)'(0) + \int_0^c (c - s)\left[K(R(s)) + h'(R(s))\varphi_r(\sigma(s))\right]ds,$$ or, equivalently, $$h(b) \leq h(a) + ch'(a)\cos \alpha + \int_0^c (c - s)\left[K(R(s)) + h'(R(s))\varphi_r(\sigma(s))\right]ds,$$ where $c = d(q, x)$ and $\alpha = \angle(\nabla r_p(q), \sigma'(0)).$ We can write this inequality in terms of the interior angle $\beta = \pi - \alpha$ at $q$ determined by $-\nabla r_p(q)$ and $\sigma'(0)$ $$h(b) \leq h(a) - ch'(a)\cos \beta + \int_0^c (c - s)\left[K(R(s)) + h'(R(s))\varphi_r(\sigma(s))\right]ds,$$
\end{theorem}

\begin{proof} Let $U(s) = u(\sigma(s))$. Observe that \begin{align*}
        U''(s) &= \Hess u(\sigma'(s), \sigma'(s))\\ &= \MHess u(\sigma'(s), \sigma'(s)) + \left<\nabla u (\sigma(s)), \nabla \varphi (\sigma(s))\right>|\sigma'(s)|^2 \\ &= \MHess u(\sigma'(s), \sigma'(s)) + h'(R(s))\varphi_r(\sigma(s))\\ &\leq K(R(s)) + h'(R(s))\varphi_r(\sigma(s)) 
    \end{align*}  for $s \in [0, c]$. Set $$\Psi(s) = K(R(s)) + h'(R(s))\varphi_r(\sigma(s)).$$ Then $$U''(s) \leq \Psi(s) \, \, \, \, \, \text{for} \, \, s \in [0, c].$$ Integrating over $[0, s]$ gives $$U'(s) - U'(0) \leq \int_0^s \Psi(t)dt,$$ from which we get $$U'(s) \leq U'(0) + \int_0^s \Psi(t)dt.$$ Integrating again over $[0, c]$ gives $$U(c) - U(0) \leq cU'(0) + \int_0^c\int_0^s \Psi(t)dtds.$$ Switching the order of integration yields the integration region $0 \leq t \leq s \leq c$, a triangle, so that $$\int_0^c\int_0^s \Psi(t)dtds = \int_0^c\int_t^c\Psi(t)dsdt = \int_0^c (c - t)\Psi(t)dt. $$ Thus $$U(c) \leq U(0) + cU'(0) + \int_0^c (c - t)\Psi(t)dt.$$ Substituting the definitions of $U$ and $\Psi$ in the preceding inequality gives $$u(x) \leq u(q) + c(u\circ\sigma)'(0) + \int_0^c (c - s)\left[K(R(s)) + h'(R(s))\varphi_r(\sigma(s))\right]ds.$$ Since $(u \circ \sigma)'(s) = du_{\sigma(s)}(\sigma'(s)) = \left<\nabla u(\sigma(s)), \sigma'(s)\right>$, at $s = 0$ we get $$(u \circ \sigma)'(0) =  \left<\nabla u(\sigma(0)), \sigma'(0)\right> = h'(r_p(q)) \left<\nabla r_p(q), \sigma'(0)\right> = h'(a)\cos \alpha$$ where $\alpha = \angle(\nabla r_p(q), \sigma'(0)).$ Hence the equivalent inequality is $$h(b) \leq h(a) + ch'(a)\cos \alpha + \int_0^c (c - s)\left[K(R(s)) + h'(R(s))\varphi_r(\sigma(s))\right]ds.$$
\end{proof}

\begin{remark} We refer to Theorem \ref{thm:3.1} as a hinge comparison due to the geometric configuration: we have a hinge angle $\beta = \pi - \alpha$ at the point $q$, where $\beta$ is the angle determined by $-\nabla r_p(q)$ and $\sigma'(0)$, the two unit directions pointing away from $q$, one along the minimizing radial geodesic toward $p$ and the other along the minimizing geodesic $\sigma$ toward $x$, respectively. The corresponding side lengths are $a = d(p, q)$, $b = d(p, x)$, and $c = d(q, x)$.
\end{remark}

\begin{proposition}[Weighted Law of Cosines]\label{prop:3.3} Assume the hypotheses of Theorem \ref{thm:3.1} and suppose that there exists a continuous function $K:(0, \infty) \rightarrow \mathbb R$ such that  $$\MHess\left(\frac{1}{2}r_p^2\right) \leq K(r_p)g$$ along $\sigma$. Define $$J_{p, K}(q, x; \sigma) = \int_0^c (c - s)K(R(s))ds\, \, \, \, \, \text{and}\, \, \, \, \, I_p(q, x; \sigma) = \int_0^c (c - s)R(s)\varphi_r(\sigma(s))ds$$ where $R(s) = r_p(\sigma(s))$ and $\varphi_r = \left<\nabla \varphi, \nabla r_p\right>$. Then the following hold.

\begin{enumerate}
    \item[(i)] (General weighted law of cosines) $$\begin{cases}b^2 \leq a^2 - 2ac\cos \beta + 2J_{p, K}(q, x; \sigma) + 2I_p(q, x; \sigma) \\ \cos\beta \leq \dfrac{a^2 - b^2 + 2J_{p,K}(q, x; \sigma) + 2I_p(q, x; \sigma)}{2ac}\end{cases}.$$

    \item[(ii)] (Weighted law of cosines under uniform upper bound on $\MHess(r_p^2/2)$) If $K(r) \leq K_0$ for some constant $K_0$ along $\sigma$, then $$\begin{cases}b^2 \leq a^2 - 2ac \cos \beta + K_0c^2 + 2I_p(q, x; \sigma)\\
    \cos \beta \leq \dfrac{a^2 - b^2 + K_0c^2 + 2I_p(q, x; \sigma)}{2ac}
    \end{cases}.$$

    \item[(iii)] (Weighted law of cosines under uniform upper bound on $\MHess(r_p^2/2)$ and radial density control) If in addition to $K(r) \leq K_0$, we have $|\nabla \varphi| \leq \dfrac{c_0}{r_p}$ for some $c_0 \geq 0$ along $\sigma$, then $$\begin{cases}b^2 \leq a^2 - 2ac\cos \beta + (K_0 + c_0)c^2 \\ \cos\beta \leq \dfrac{a^2 - b^2 + (K_0 + c_0)c^2}{2ac}\end{cases}.$$
\end{enumerate}

\end{proposition}

\begin{proof}
Note that with the radial potential $h(r) = r^2/2$, $h'(R(s)) = R(s)$. So (i) follows by applying Theorem \ref{thm:3.1} with $h(r_p) = r_p^2/2$ and multiplying the resulting inequality by $2$. Part (ii) follows from $$J_{p, K}(q, x; \sigma) \leq K_0 \int_0^c (c - s)ds = \dfrac{K_0}{2}c^2.$$ For (iii), we use the bound on $J_{p, K}(q, x; \sigma)$ obtained above along with \begin{align*}|I_p(q, x; \sigma)| &\leq \int_0^c (c - s)R(s)|\varphi_r(\sigma(s))|ds \\ &= \int_0^c (c - s)R(s)|\left<\nabla \varphi(\sigma(s)), \nabla r_p(\sigma(s))\right>|ds \\ &\leq \int_0^c (c - s)R(s)|\nabla \varphi(\sigma(s))||\nabla r_p(\sigma(s))|ds \\ &\leq \int_0^c (c - s)r_p(\sigma(s)) \dfrac{c_0}{r_p(\sigma(s))}ds \\ &\leq c_0\int_0^c (c - s)ds \\ &= c_0\dfrac{c^2}{2}.\end{align*}
\end{proof}

\begin{remark}\label{rm:3.4}
Since $p\ne q$ and $q\ne x$, we have $a,c>0$, so the
angle estimates are well-defined. The degenerate triangle case $p\rightarrow q\rightarrow x$ is still allowed; in this case, $a+c=b$ and $\beta=\pi$. Thus, our assumptions exclude zero-length sides but do not exclude degenerate triangles. Moreover, we state the angle estimates in terms of the interior angle $\beta$, as it is the natural hinge angle at $q$. The corresponding estimates for the exterior angle $\alpha$ at $q$ follow immediately from $\alpha=\pi-\beta$.

We also note that our hinge comparison is somewhat different from the classical hinge comparison in Toponogov's theorem \cite{Meyer}. In our setting, we fix the basepoint $p$ and consider the radial distance function $r_p$ along the geodesic $\sigma$ from $q$ to $x$. Thus, rather than comparing a hinge with a corresponding hinge in a model space, our comparison comes from an upper bound on the modified Hessian of a radial function based at $p$.
\end{remark}

\begin{remark}
Note that the inequalities of Proposition \ref{prop:3.3} (iii) hold under the hypotheses $\overline{\mathrm{sec}}_\varphi \geq 0$, $a_0 \leq \varphi \leq 0$, and $|\nabla \varphi| \leq c_0/r_p$, as these hypotheses provide a uniform bound $K_0$ on the modified Hessian of the squared distance function (see \cite[Theorem 3.1]{NgPaper1} for details). Furthermore, one can show that under the same nonnegative weighted sectional curvature and corresponding density controls that $$\MHess(r_p^2/2) \leq K\tilde{g},$$ where the conformal metric $\tilde{g} = e^{-2\varphi}g$. This leads us to consider a ``mixed geodesic triangle,'' i.e., its sides from $p$ to $q$ and from $p$ to $x$ are minimizing $g$-geodesics, while its side from $q$ to $x$ is a minimizing $\tilde{g}$-geodesic. This provides us a ``Warped Law of Cosines'': $$b^2\leq a^2+k_1\widetilde c^{\,2}
-2k_2a\widetilde c\cos\beta,$$
where $\tilde{\sigma}:[0, \tilde{c}] \rightarrow M$ is a unit-speed minimizing geodesic with respect to $\tilde{g}$ from $q$ to $x$ and $$\tilde{c} = L(\tilde{g}(\tilde{\sigma})), \qquad \beta=\angle\bigl(-\nabla r_p(q),\widetilde\sigma'(0)\bigr), \qquad k_1=(1+C_0)e^{-2a_0}, \qquad k_2=
\begin{cases}
e^{-a_0},&0\leq\beta<\dfrac{\pi}{2},\\[2mm]
e^{2a_0},&\dfrac{\pi}{2}\leq\beta\leq\pi
\end{cases}.$$ We hope that these estimates may eventually provide a framework for establishing Toponogov-type comparison results in this setting.
\end{remark}

\section{Lower Excess Estimates from the Weighted Law of Cosines}\label{sec:4}

We now establish some lower excess estimates in the setting of Theorem \ref{thm:3.1}. Recall that in the setup of Section \ref{sec:3}, $(M^n, g, \varphi)$ is a Riemannian manifold with density, $p \in M$ is a fixed basepoint of the distance function $r = r_p$, and $$a = r_p(q) = d(p, q), \, \, \, \, \, b = r_p(x) = d(p, x), \, \, \, \, \, \text{and} \, \, \, \, \, c = d(q, x) = L(\sigma),$$ where $\sigma$ is a unit-speed minimizing geodesic from $q$ to $x$. We consider the excess at $q \in M$ given by $$e(q) = d(p, q) + d(q, x) - d(p, x) = a + c - b.$$ The excess $e(q)$ measures the excess length of the path from $p$ to $q$, then from $q$ to $x$, relative to the direct distance from $p$ to $x$.  Note that $e(q) \geq 0$ by the triangle inequality. We allow for the degenerate case $a + c = b$ with interior hinge angle $\beta = \pi$ at $q$. For complete manifolds, $e(q) = 0$ if and only if $q$ belongs to a minimizing geodesic from $p$ to $x$. We will often write $e = a + c - b$ for brevity. 

The following lower excess estimates require a radial potential bound $K(r_p)$ along $\sigma$ on the modified Hessian of the modified squared distance function $u = h \circ r = r_p^2/2$. 

\begin{proposition}[Excess Inequality]\label{prop:4.1}
    Assume the hypotheses of Theorem \ref{thm:3.1}, and suppose that there exists a continuous function $K:(0, \infty) \rightarrow \mathbb R$  such that $\MHess\left(r_p^2/2\right) \leq K(r_p)g$ along $\sigma$. Then the excess $e = a + c - b$ satisfies \begin{equation}\label{eq:4.1}2(a + c)e - e^2 \geq 2ac(1 + \cos\beta) + c^2 - 2J_{p, K}(q, x; \sigma) - 2I_p(q, x; \sigma),\end{equation} where $J_{p, K}(q, x; \sigma)$ and $I_p(q, x; \sigma)$ are the integrals defined in Proposition \ref{prop:3.3}. 
\end{proposition}

\begin{proof}
    From the definition of excess at $q$, we have $b = a + c - e,$ and squaring both sides gives $$b^2 = (a + c)^2 - 2(a + c)e + e^2.$$ Since $u = h(r) = r^2/2$, by Proposition \ref{prop:3.3}(i) we also have $$b^2 \leq a^2 - 2ac\cos\beta + 2J_{p, K}(q, x; \sigma) + 2I_p(q, x; \sigma).$$ Putting these together yields $$(a + c)^2 - 2(a + c)e + e^2 \leq a^2 - 2ac\cos\beta + 2J_{p, K}(q, x; \sigma) + 2I_p(q, x; \sigma).$$ Solving for $2(a + c)e - e^2$ gives inequality \eqref{eq:4.1}.   
\end{proof}

\begin{corollary}[Lower Excess Estimates]
    Under the hypotheses of Proposition \ref{prop:4.1}, the following estimates hold. 

    \begin{enumerate}
        \item[(a)] (General lower excess estimates) $$\begin{cases}
            e \geq \dfrac{ac}{a + c}(1 + \cos\beta) + \dfrac{c^2}{2(a + c)} - \dfrac{J_{p, K}(q, x; \sigma)}{a + c} - \dfrac{I_p(q, x; \sigma)}{a + c} & (1) \\
            e \geq a + c - \sqrt{(a + c)^2 - \left(2ac(1 + \cos\beta) + c^2 - 2J_{p, K}(q, x; \sigma) - 2I_p(q, x; \sigma)\right)} & (2)\\
            1 + \cos\beta \leq \dfrac{2(a + c)e - e^2 - c^2 + 2J_{p, K}(q, x; \sigma) + 2I_p(q, x; \sigma)}{2ac} & (3)
        \end{cases}.$$

        \item[(b)] (Lower excess estimates under uniform bound on $K(r)$) 
        
        Suppose, in addition, that $K(r) \leq K_0$ for some constant $K_0$ along $\sigma$. Then \begin{equation}\label{eq:4.2} 2(a + c)e - e^2 \geq 2ac(1 + \cos\beta) - (K_0 - 1)c^2 - 2I_p(q, x; \sigma),\end{equation} and so $$\begin{cases}
            e \geq \dfrac{ac}{a + c}(1 + \cos\beta) - \dfrac{(K_0 - 1)c^2}{2(a + c)}  - \dfrac{I_p(q, x; \sigma)}{a + c} & (1)\\
            e \geq a + c - \sqrt{(a + c)^2 - \left(2ac(1 + \cos\beta) -(K_0 -1)c^2 - 2I_p(q, x; \sigma)\right)} & (2)\\
            1 + \cos\beta \leq \dfrac{2(a + c)e - e^2 + (K_0 - 1)c^2 + 2I_p(q, x; \sigma)}{2ac} & (3)
        \end{cases}.$$
        \item[(c)] (Lower excess estimates under uniform radial bound and additional density control) 

If in addition to the hypothesis of (b) we have $|\nabla \varphi| \leq c_0/r_p$ for some constant $c_0$ along $\sigma$, then \begin{equation}\label{eq:4.3}2(a + c)e - e^2 \geq 2ac(1 + \cos\beta) - (K_0 -1 + c_0)c^2,\end{equation} and so $$\begin{cases}
            e \geq \dfrac{ac}{a + c}(1 + \cos\beta) - \dfrac{K_0 -1 + c_0}{2(a + c)}c^2 & (1)\\
            e \geq a + c - \sqrt{(a + c)^2 - \left(2ac(1 + \cos\beta) -(K_0 -1 + c_0)c^2\right)} & (2)\\
            1 + \cos\beta \leq \dfrac{2(a + c)e - e^2 + (K_0 - 1 + c_0)c^2}{2ac} & (3)
        \end{cases}.$$ 
    \end{enumerate}

    \noindent The lower bound in (a)(2) is positive whenever the right-hand side of \eqref{eq:4.1} is positive. If the right-hand side of \eqref{eq:4.1} is $0$, then (a)(2) reduces to the trivial estimate $e \geq 0$. Similarly, the lower bound of (b)(2) is positive whenever the right-hand side of \eqref{eq:4.2} is positive, and the analogous statement for (c)(2) follows from \eqref{eq:4.3}.
\end{corollary}

\begin{proof}
Estimate (a)(1) follows from $2(a + c)e - e^2 \leq 2(a + c)e$ as $e\geq 0$. The stronger estimate (a)(2) follows from solving for $e$ directly in \eqref{eq:4.1} by completing the square on the left-hand side of \eqref{eq:4.1}. The angle estimate (a)(3) follows from Remark \ref{rm:3.4} and rearranging \eqref{eq:4.1}.
\end{proof}

\begin{example}[Recovery of Euclidean Law of Cosines]
  In the case of constant density $\varphi = k$ in Euclidean space $M^n = \mathbb R^n$, we have $$\MHess(r_p^2/2) = g, \qquad I_p(q,x) = 0.$$ Since $K(r) = 1$, we obtain $$J_{p, K}(q, x; \sigma) = \int_0^c(c -s)ds = \dfrac{c^2}{2}.$$ Since the modified Hessian bound of $1$ is attained, the comparison in Proposition \ref{prop:3.3}(i) is an equality, and we recover the Euclidean Law of Cosines $$b^2 = a^2 - 2ac\cos\beta + c^2.$$ In this case, the corresponding excess estimate becomes the Euclidean excess identity $$2(a + c)e - e^2 = 2ac(1 + \cos \beta).$$
\end{example}

\section{Abresch-Gromoll Excess Estimate under Nonnegative Weighted Sectional Curvature}\label{sec:5}

Given $p, q \in M$, we establish an upper bound on the excess at $x \in M$, $$e(x) = d(p, x) + d(x, q) - d(p, q)$$ under nonnegative weighted sectional curvature and density control. The proof makes use of the weighted radial Laplacian comparison with a radial barrier argument to obtain this upper bound.

\begin{theorem}\label{thm:5.1}
    Let $(M^n, g, \varphi)$ be a complete $n$-dimensional Riemannian manifold with density with $n \geq 3$. Suppose $$\overline{\mathrm{sec}}_\varphi \geq 0, \qquad a_0 \leq \varphi \leq 0, \qquad |\nabla \varphi| \leq \frac{c_0}{r_p},$$ where $r_p(x) = d(p, x)$. Define $$\alpha = (n - 1)e^{-2a_0}, \qquad \beta = 2c_0, \qquad A = \alpha + \beta.$$ Let $\gamma:[0, L] \rightarrow M$ be a unit-speed minimizing geodesic from $p$ to $q$, and define $$h = d(x, \gamma), \qquad s = \min \{d(p, x), d(q, x)\}.$$ If $h \leq \dfrac{s}{2}$, then \begin{equation}\label{eq:5.1} e(x) \leq \mathcal C_A\left(\dfrac{h^{A + 1}}{s}\right)^{\frac{1}{A}}, \qquad \qquad \mathcal C_A = \dfrac{2A}{A - 1}\left(\dfrac{2A}{A + 1}\right)^{\frac{1}{A}}.\end{equation}
\end{theorem}

\begin{remark}\label{rm:5.2}
If $h = 0$, then $x \in \gamma$, and so $e(x) = 0$. Thus, the estimate is immediate in this case. Moreover, the restriction $h \leq s/2$ is needed in the proof and describes the case we are most interested in, namely, the case where triangles are close to being degenerate (i.e., thin triangles). Indeed, the height $h$ is small relative to the minimum of the side lengths $s$, and the crude estimate $e(x) \leq 2h$, which always holds by considering $z \in \gamma$ such that $d(x, z) = h$, shows that the excess is controlled by the height and hence by $s$. On the other hand, if $h > s/2$, then $$\mathcal C_A\left(\frac{h^{A+1}}{s}\right)^{1/A} = \mathcal C_A h\left(\frac {h}{s}\right)^{1/A} > \mathcal C_A2^{-1/A}h \geq 2h \geq e(x),$$ where the second to last inequality follows from
$$\left(1-\frac1A\right) \leq \left(\dfrac{1}{e}\right)^{1/A} < \left(\frac{1}{2}\right)^{1/A} <\left(\frac{A}{A+1}\right)^{1/A}, \qquad A>1.$$ Thus, the same estimate \eqref{eq:5.1} holds also when $h > s/2$. 
\end{remark}

\begin{remark}
The assumption $n\geq 3$ guarantees that $A\geq 2$ under our
hypotheses. In particular, the parameter $$A=(n-1)e^{-2a_0}+2c_0, \qquad a_0 \leq 0, \, \, c_0 \geq 0$$
depends on the dimension as well as the bounds on the density and
its gradient. We note that the estimate \eqref{eq:5.1} holds for the more general range $A>1$. This follows from the barrier
argument in the proof of Theorem \ref{thm:5.1} when $h\leq s/2$,
and from the algebraic argument in Remark \ref{rm:5.2} above when $h>s/2$.

In the zero-density case $\varphi\equiv 0$, we have $a_0=c_0=0$,
and hence $A=n-1$. Our estimate \eqref{eq:5.1} becomes
$$
e(x)\leq \mathcal C_{n-1}
\left(\frac{h^n}{s}\right)^{\frac{1}{n-1}}, \qquad \text{where}\qquad 
\mathcal C_{n-1}
=
\frac{2(n-1)}{n-2}
\left(\frac{2(n-1)}{n}\right)^{\frac{1}{n-1}}.
$$
For $n=3$, this gives
$$
e(x)\leq \frac{8}{\sqrt{3}}
\left(\frac{h^3}{s}\right)^{\frac12},
$$
which improves the constant $8$ in the classical
Abresch--Gromoll estimate
$$
e(x)\leq 8\left(\frac{h^3}{s}\right)^{\frac12}.
$$ but with a stronger curvature assumption than the Ricci curvature assumption
in the classical result, as zero-density nonnegative weighted sectional curvature is equivalent to $\sec \geq 0$. More generally, our estimate has the same exponent and scaling
as the classical Abresch--Gromoll estimate, but with a different
constant. Recall that the classical estimate holds under
$n\geq 3$ and $\operatorname{Ric}\geq 0$
(see \cite[Theorem 2.5.4]{DaiWei} and \cite{AbreschGromoll}).

Wei--Wylie also established an excess estimate under
$\operatorname{Ric}_f\geq 0$ and $|f|\leq k$
\cite[Theorem 7.1]{WeiWylie09}. Their estimate states that
$$e(x)\leq 2\left(\frac{n+4k-1}{n+4k-2}\right)
\left(\dfrac{1}{2}Ch^{n + 4k}\right)^{\frac{1}{n+4k-1}}.$$ where $$
C = 2\left(\dfrac{n + 4k -1}{n + 4k}\right)\left(\dfrac{1}{d(p, x) - h(x)} + \dfrac{1}{d(q,x) - h(x)}\right), \qquad h(x) < \min\{d(p, x), d(x, q)\}.$$
Thus, as in our estimate, the exponent in this setting depends on the dimension
and bounds on the density. In the Wei--Wylie estimate, the exponent $n+4k$ corresponds to our exponent $A + 1$, so that $A = n + 4k -1$.
The two estimates therefore have a similar form, but arise
from different curvature assumptions and comparison methods.
The Wei--Wylie estimate uses a lower bound on the
Bakry--Emery Ricci tensor and a bound on the potential, while
ours uses a modified Hessian comparison under
nonnegative weighted sectional curvature and control of the
density and its gradient.
\end{remark}

\begin{lemma}[Weighted radial Laplacian comparison]\label{lem:5.4} Assume the hypotheses of Theorem \ref{thm:5.1}. For any $y \in M$ and  wherever $r_y$ and $r_p > 0$, \begin{equation}\label{eq:5.2}\Delta_\mu r_y \leq \dfrac{\alpha}{r_y} + \dfrac{\beta}{r_p},\end{equation} where $\alpha = (n - 1)e^{-2a_0}$ and $\beta = 2c_0.$
\end{lemma}

\begin{proof}[Proof of Lemma \ref{lem:5.4}]
Recall that the tangential part of the modified Hessian satisfies $$\Hess r_y - \varphi_{r_y}g \leq \dfrac{e^{-2\varphi}}{r_y}g$$ (see the proof of \cite[Theorem 3.1]{NgPaper1}. Taking the trace over an orthonormal basis $\{E_i\}_{i=1}^n$ where $E_n = \nabla r_y$ together with $a_0 \leq \varphi \leq 0$ gives $$\Delta r_y = \tr_g\Hess r_y \leq (n - 1) \dfrac{e^{-2a_0}}{r_y} + (n-1)\varphi_{r_y}.$$ The weighted Laplacian with respect to smooth measure $\mu$ is then \begin{align*}
    \Delta_\mu r_y = \Delta r_y -(n + 1)\left<\nabla \varphi, \nabla r_y\right> &= \Delta r_y - (n + 1)\varphi_{r_y} \\ &\leq \left(\frac{(n -1)e^{-2a_0}}{r_y} + (n-1)\varphi_{r_y}\right) -(n + 1)\varphi_{r_y} \\ &= \frac{(n -1)e^{-2a_0}}{r_y} - 2\varphi_{r_y} \\ &\leq \frac{(n -1)e^{-2a_0}}{r_y} + 2\dfrac{c_0}{r_p} = \dfrac{\alpha}{r_y} + \dfrac{\beta}{r_p}
\end{align*} as desired. 
\end{proof}

To obtain the weighted Laplacian comparison at points in the $\mathrm{Cut}(y)$, we
use Calabi support functions. Fix a point $z\neq y$ and let
$\gamma:[0,\ell]\to M$ be a minimizing unit-speed geodesic from $y$ to
$z$. For sufficiently small $\epsilon>0$, set
$y_\epsilon=\gamma(\epsilon)$ and define $$\rho_{y,\epsilon}(\cdot)
=d(y,y_\epsilon)+d(y_\epsilon,\cdot).$$ Then $\rho_{y,\epsilon}\geq r_y$, with equality at $z$. Moreover, $\rho_{y,\epsilon}$ is smooth in a neighborhood of $z$ for sufficiently small $\epsilon>0$, and the weighted Laplacian comparison of $\rho_{y, \epsilon}$ is \begin{equation}\label{eq:5.3} \Delta_\mu \rho_{y, \epsilon}(\cdot) = \Delta_\mu (d(y, y_\epsilon) + d(y_\epsilon, \cdot)) = \Delta_\mu r_{y_\epsilon}(\cdot) \leq \dfrac{\alpha}{r_{y_\epsilon}(\cdot)} + \dfrac{\beta}{r_p(\cdot)} = \dfrac{\alpha}{\rho_{y, \epsilon}(\cdot) - \epsilon} + \dfrac{\beta}{r_p(\cdot)}\end{equation} on this neighborhood. Although the right-hand side of \eqref{eq:5.2} is singular at $z=p$, in our arguments we will work in regions where $r_p$ is bounded away from zero, so that the right-hand side remains finite. In other words, the neighborhood of $z$ on which \eqref{eq:5.3} holds is chosen sufficiently small to avoid $p$. Thus, Calabi support functions are only needed to handle the cut locus of distance functions $r_y$.\\

\noindent \textit{Proof Strategy for Theorem \ref{thm:5.1}}
\begin{enumerate}
    \item We first obtain bounds for the weighted Laplacians
    $\Delta_\mu e$ and $\Delta_\mu r_x$ on a ball centered at $x$, where $e$ is the excess function and $r_x=d(x,\cdot)$.
    We choose the ball so that it stays away from $p$ and $q$,
    allowing us to use these bounds throughout the ball.

    \item Using these bounds, we construct a function
    $G_\delta=g_\delta(r_x)$ that depends only on the distance
    from $x$. We choose $g_\delta$ so that the weighted Laplacian $\Delta_\mu G_\delta$ is strictly larger than the weighted Laplacian $\Delta_\mu e$. This function will be used to compare with the excess.

    \item We then argue by contradiction on an annulus $\Omega_{t, R} = \left\{t \leq r_x(y) \leq R\right\}$ centered at $x$. If the excess $e(x)$ is too large, i.e., if $e(x) > 2t + g_\delta(t),$ the difference $G_\delta-e$ is positive somewhere in the annulus but nonpositive on its boundary. It must therefore attain a positive maximum at an interior point.

    \item We show that this interior maximum is impossible. Since the distance functions involved may not be smooth at
    that point, we replace them locally by smooth supports that agree with them at the maximum. The Laplacian bounds
    then give a contradiction to the existence of a positive
    interior maximum of $G_\delta-e$.

    \item Finally, we let the small perturbation $\delta$ tend
    to zero, let the outer radius approach the height $h$, and
    optimize over $t$. This gives the desired upper bound for the excess.
\end{enumerate}

\begin{proof}[Proof of Theorem \ref{thm:5.1}]

Assume the hypotheses of the theorem. We first show that the weighted Laplacian of both the excess function at $x$ and the distance function centered $x$ can be bounded uniformly on a ball of radius $R$ centered at $x$, $B(x, R)$. The excess is $$e = r_p + r_q - L.$$ By linearity of $\Delta_\mu$ and Lemma \ref{lem:5.4}, we have  \begin{align*}
    \Delta_\mu e = \Delta_\mu r_p + \Delta_\mu r_q  &\leq \left(\dfrac{\alpha}{r_p} + \dfrac{\beta}{r_p}\right) + \left(\dfrac{\alpha}{r_q} + \dfrac{\beta}{r_p}\right) \\ &= \dfrac{\alpha + 2\beta}{r_p} + \dfrac{\alpha}{r_q}. 
\end{align*} Since $h < s$, we choose the radius $R$ of $B(x, R)$ where $h < R < s$. If $y \in B(x, R)$, then $$r_p(y) = d(p, y) \geq d(p, x) - d(x, y) > s - R.$$ Similarly, $$r_q(y) = d(q, y) > s - R.$$ Thus $$\Delta_\mu e \leq \dfrac{\alpha + 2\beta}{r_p} + \dfrac{\alpha}{r_q} < \dfrac{\alpha + 2\beta}{s - R} + \dfrac{\alpha}{s - R} = \dfrac{2(\alpha + \beta)}{s - R} = \dfrac{2A}{s - R} \qquad \text{on} \, \, \, \, \, B(x, R).$$ In a similar manner, $$\Delta_\mu r_x \leq \dfrac{\alpha}{r_x} + \dfrac{\beta}{r_p} \leq \dfrac{\alpha}{r_x} + \dfrac{R}{r_x}\dfrac{\beta}{s - R} = \left(\alpha + \dfrac{R\beta}{s - R}\right)\dfrac{1}{r_x}\qquad \text{on} \, \, \, \, \,   B(x, R).$$ Define $$\lambda_R = \max\left\{A, \alpha + \dfrac{\beta R}{s - R}\right\} \qquad \text{and} \qquad B_R = \max \left\{\dfrac{2A}{s - R}, \dfrac{4A}{s}\right\}.$$ Since $h \leq s/2$, we have $$ \dfrac{h}{s -h} \leq 1 \implies \lim_{R \downarrow h} \lambda_R = A \qquad \text{and} \qquad \dfrac{1}{s - h} \leq \dfrac{2}{s} \implies \lim_{R \downarrow h} B_R = \dfrac{4A}{s}.$$ Putting these together, we have $$\Delta_\mu r_x \leq \dfrac{\lambda_R}{r_x} \qquad \text{and} \qquad \Delta_\mu e \leq B_R,$$ and in the limit $R \downarrow h$, $$\Delta_\mu r_x \leq A \qquad \text{and} \qquad \Delta_\mu e \leq \dfrac{4A}{s}.$$

We now construct a radial barrier function $$G(y) = g(r_x(y))$$ where $g:[0, R] \rightarrow \mathbb R$ so that $$\Delta_\mu(G - e) \geq 0.$$ To get this, we want $\Delta_\mu G$ to be bounded below by $B_R$. Write $r = r_x$. By the chain rule, $$\Delta_\mu G(y) = \Delta_\mu (g(r)) = g''(r) + \Delta_\mu rg'(r).$$ Since $\Delta_\mu r \leq \lambda_R/r$ on $B(x, R)$, we want $g'(r) < 0$ so that $$\Delta_\mu G \geq g''(r) + \dfrac{\lambda_R}{r}g'(r).$$ So to obtain $g$ with the desired properties, we solve the ordinary differential equation $$\begin{cases}
    g''(u) + \dfrac{\lambda_R}{u}g'(u) = B_R \\
    g(R) = 0\\
    g'(R) = 0
\end{cases}$$ where the boundary conditions are chosen such that the barrier vanishes at the outer radius $R$ and to fix an integration constant. But rather than solving this ODE, we solve a perturbed ODE where the right-hand side is $B_R + \delta$, where $\delta > 0$. This gives a barrier function $G_\delta = g_\delta(r)$ satisfying the strict inequality $$\Delta_\mu (G_\delta - e) \geq \delta > 0,$$ allowing us to rule out an interior maximum of $G_\delta - e$ by the second derivative test. We fix $\delta > 0$ and solve $$\begin{cases}
    g_\delta''(u) + \dfrac{\lambda_R}{u}g_\delta'(u) = B_R + \delta & (*) \\
    g_\delta(R) = 0\\
    g_\delta'(R) = 0.
\end{cases}$$ Multiplying $(*)$ by $u^{\lambda_R}$, we can write it as $$\dfrac{d}{du}(u^{\lambda_R}g_\delta'(u)) = u^{\lambda_R}(B_R + \delta).$$ Integrating this gives $$g_\delta'(u) = \dfrac{B_R + \delta}{\lambda_R + 1}\left(u - \dfrac{R^{\lambda_R + 1}}{u^{\lambda_R}}\right),$$ which is strictly negative on $0 < u < R$ as desired. In other words, $g_\delta$ is strictly decreasing as $u$ increases toward $R$. But since $g_\delta(R) = 0$, it follows that $g_\delta(u) > 0$ on $0 < u < R$. Integrating again gives $$g_\delta(u) = \dfrac{B_R + \delta}{2(\lambda_R + 1)}\left(u^2 + \dfrac{2R^{\lambda_R + 1}}{\lambda_R -1}u^{1 - \lambda_R} - \dfrac{\lambda_R + 1}{\lambda_R -1}R^2\right).$$ 

We now show that $$e(x) \leq 2t + g_\delta(t), \qquad 0 < t < h < R.$$ Suppose for contradiction that there exists $t \in (0, h)$ and $R \in (h, s)$ such that $$e(x) > 2t + g_\delta(t).$$ Consider the closed and bounded annulus centered at $x$ $$\Omega_{t, R} = \{y \in M \mid t \leq r_x(y) \leq R\}.$$ Since $M$ is complete, it follows by Hopf-Rinow that $\Omega_{t, R}$ is compact. Define $$V_\delta(y) = G_\delta(y) - e(y).$$ Since this function is continuous away from $x$, and $\Omega_{t, R}$ is compact, $V_\delta$ achieves its maximum on $\Omega_{t, R}$. We claim that $\max_{\Omega_{t, R}} V_\delta$ does not occur on $\partial \Omega_{t, R}$. For the inner boundary, suppose $y \in M$ such that $r_x(y) = d(x, y) = t$. Since the excess function is $2$-Lipschitz, we have $$|e(x) - e(y)| \leq 2d(x, y) = 2t \implies e(x) - 2t \leq e(y).$$ This together with $e(x) > 2t + g_\delta(t)$ gives $$e(y) \geq e(x) - 2t > (2t + g_\delta(t)) - 2t = g_\delta(t),$$ i.e., $e(y) > g_\delta(r_x(y))$. Equivalently, $V_\delta < 0$ on the inner boundary $r_x(y) = t$. For the outer boundary, we have $r_x(y) = R$ and $g_\delta(R) = 0$. So $$V_\delta(y) = g_\delta(r_x(y)) - e(y) = g_\delta(R) - e(y) = -e(y) \leq 0.$$ Thus $V_\delta \leq 0$ on the outer boundary $r_x(y) = R$, and, more generally, on $\partial \Omega_{t, R}.$

We now show there exists an interior point of $\Omega_{t, R}$ where $V_\delta > 0$. Let $z \in \gamma$ such that $d(x, z) = h$. Since $t < h < R$, we have $z \in \mathrm{Int}\Omega_{t, R}$. Since $e(z) = 0$, $g_\delta(r_x(z)) = g_\delta(h) > 0$ by our observations above with $0 < h < R$. Therefore $$V_\delta(z) = g_\delta(h) - e(z) > 0, \qquad z \in \mathrm{Int}\Omega_{t, R}.$$ Since $V_\delta \leq 0$ on $\partial \Omega_{t, R}$ and $V_\delta > 0$ on some $z \in \mathrm{Int}\Omega_{t, R}$, it follows that $\max_{\Omega_{t, R}} V_\delta > 0$ at an interior point $y_*$. So $t < r_x(y_*) < R$.

The goal now is to derive a contradiction at this maximum point $y_* = \max_{\Omega_{t, R}} V_\delta$. Note \begin{align*} V_\delta(y_*) &= g_\delta(r_x(y_*)) - e(y_*) \\ &= g_\delta(r_x(y_*)) - \left(r_p(y_*) + r_q(y_*) - L\right).\end{align*} We need to account for the case where $V_\delta$ may fail to be smooth at $y_*$ as it may be a cut point of the distance functions $r_x, r_p,$ or $r_q$. To handle this case, we construct a smooth lower support $V_\epsilon$ for $V_\delta$ via smooth upper supports for the three distance functions. We choose unit-speed minimizing geodesics from $x, p,$ and $q$ to $y_*$. For $\epsilon_x, \epsilon_p, \epsilon_q > 0$ sufficiently small, let $x_{\epsilon_x}, p_{\epsilon_p},$ and $q_{\epsilon_q}$ be points a distance $\epsilon_x, \epsilon_p,$ and $\epsilon_q$, respectively, along the corresponding minimizing geodesics. Then define supports $$\rho_x(y) = d(x, x_{\epsilon_x}) + d(x_{\epsilon_x}, y) = \epsilon_x + d(x_{\epsilon_x}, y),$$ with $\rho_p(y)$ and $\rho_q(y)$ defined similarly. Note that these are shifted distance functions along their corresponding geodesics. Define $$e_\epsilon = \rho_p + \rho_q - L \qquad \text{and} \qquad V_\epsilon = g_\delta \circ \rho_x - e_\epsilon.$$ Note $e_\epsilon \geq e$ as $\rho_p \geq r_p$ and $\rho_q \geq r_q$. Also note that $\rho_x(y_*) = r_x(y_*) < R$ since $y_* \in \mathrm{Int}\Omega_{t, R}$. By continuity of $\rho_x$ and $r_x$, we can restrict to a sufficiently small neighborhood of $y_*$ such that $\rho_x < R$, i.e., the neighborhood is away from the outer boundary $r_x = R$, and such that $\rho_x$ and $r_x > 0$. Then $g_\delta(\rho_x) \leq g_\delta(r_x)$ on this neighborhood as $g_\delta$ is decreasing and $\rho_x \geq r_x$. We get $$V_\epsilon = g_\delta(\rho_x) - e_\epsilon  \leq g_\delta(r_x) - e = V_\delta$$ on a neighborhood of $y_*$ and $$V_\epsilon(y_*) = g_\delta(\rho_x(y_*)) - (\rho_p(y_*) + \rho_q(y_*) - L) = g_\delta(r_x(y_*)) - (r_p(y_*) + r_q(y_*) - L) = V(y_*).$$ Thus $V_\epsilon$ is a smooth lower support for $V_\delta$ at $y_*$. Since $y_*$ is a local maximum of  the $V_\delta$, it is also a local maximum of the smooth function $V_\epsilon$. 

We show that $\Delta_\mu V_\epsilon(y_*) > 0$, which will contradict $V_\epsilon$ having a smooth local maximum at $y_*$. We must compute $\Delta_\mu e_\epsilon(y_*)$ and $\Delta_\mu(g_\delta(\rho_x(y_*)))$. By Lemma \ref{lem:5.4} and the definition of $\rho_p$, \begin{align*}\Delta_\mu \rho_p(y_*) &= \Delta_\mu(\epsilon_p + d(p_{\epsilon_p}, y_*)) \\ &= \Delta_\mu r_{p_{\epsilon_p}}(y_*) \\ &\leq \dfrac{\alpha}{r_{p_{\epsilon_p}}(y_*)} + \dfrac{\beta}{r_p(y_*)} \\ &\leq \dfrac{\alpha}{r_p(y_*) - \epsilon_p} + \dfrac{\beta}{r_p(y_*)} \\ &= \alpha\left(\dfrac{1}{r_p(y_*)} + \dfrac{\epsilon_p}{r_p(y_*)(r_p(y_*) - \epsilon_p)}\right) + \dfrac{\beta}{r_p(y_*)} \\ &\leq \alpha\left(\dfrac{1}{r_p(y_*)} + \dfrac{2\epsilon_p}{(s - R)^2}\right) + \dfrac{\beta}{r_p(y_*)}
\end{align*} where the last equality follows from $$r_p(y_*) = d(p, y_*) \geq d(p, x) - d(x, y_*) \geq s - d(x, y_*) > s - R$$ and choosing $\epsilon_p$ such that $0 <\epsilon_p < \dfrac{s - R}{2}$. Similarly, for $0 < \epsilon_q < \dfrac{s- R}{2}$ $$\Delta_\mu \rho_q(y_*) \leq \alpha\left(\dfrac{1}{r_q(y_*)} + \dfrac{2\epsilon_q}{(s - R)^2}\right) + \dfrac{\beta}{r_q(y_*)}.$$ Then \begin{align*}\Delta_\mu e_\epsilon(y_*) &= \Delta_\mu \rho_p(y_*) + \Delta_\mu \rho_q(y_*) \\ &\leq \dfrac{\alpha + 2\beta}{r_p(y_*)} + \dfrac{\alpha}{r_q(y_*)} + \dfrac{2\alpha}{(s - R)^2}(\epsilon_p + \epsilon_q) \\ &\leq \dfrac{2(\alpha + \beta)}{s - R} + \dfrac{2\alpha}{(s - R)^2}(\epsilon_p + \epsilon_q) \\ &= \dfrac{2A}{s - R} + \dfrac{2\alpha}{(s - R)^2}(\epsilon_p + \epsilon_q) \\ &\leq B_R + \dfrac{2\alpha}{(s - R)^2}(\epsilon_p + \epsilon_q).\end{align*}

Again by Lemma \ref{lem:5.4} and the definition of $\rho_x$ we have $$\Delta_\mu \rho_x(y_*) \leq \dfrac{\alpha}{r_x(y_*) - \epsilon_x} + \dfrac{\beta}{r_p(y_*)}.$$ Choose $\epsilon_x$ such that $0 < \epsilon_x < t/2$, where $0 < t < h$ and $h < R < s$. Note $r_x(y_*) > t$ for all $y_* \in \mathrm{Int}\Omega_{t, R}$. We get $$r_x(y_*) - \epsilon_x > t - t/2 = t/2.$$ We also have $$r_p(y_*) > s - R \implies \dfrac{\beta}{r_p(y_*)} \leq \dfrac{\beta}{s - R} \leq \dfrac{\beta R}{(s - R)r_x(y_*)}.$$ Putting these observations together, we get \begin{align*}\Delta_\mu \rho_x(y_*) &\leq \dfrac{\alpha}{r_x(y_*) - \epsilon_x} + \dfrac{\beta}{r_p(y_*)} \\ &\leq \alpha\left(\dfrac{1}{r_x(y_*)} + \dfrac{\epsilon_x}{r_x(y_*)(r_x(y_*) - \epsilon_x)}\right) + \dfrac{\beta R}{(s - R)r_x(y_*)} \\ &\leq \dfrac{\alpha + \dfrac{\beta R}{s - R}}{r_x(y_*)} + \dfrac{2\alpha \epsilon_x}{t^2} \\ &\leq \dfrac{\lambda_R}{r_x(y_*)} + \dfrac{2\alpha \epsilon_x}{t^2}.\end{align*} It follows that \begin{align*}\Delta_\mu (g_\delta(\rho_x))(y_*) &= g_\delta''(r_x(y_*))\underbrace{|\nabla \rho_x(y_*)|^2}_{=1} + \underbrace{g_\delta'(r_x(y_*))}_{<0}\Delta_\mu \rho_x(y_*) \\ &\geq g_\delta''(r_x(y_*)) + \dfrac{\lambda_R}{r_x(y_*)}g_\delta'(r_x(y_*)) + \dfrac{2\alpha \epsilon_x}{t^2}g_\delta'(r_x(y_*)) \\ &= B_R + \delta - \dfrac{2\alpha}{t^2}|g_\delta'(r_x(y_*))|\epsilon_x.\end{align*} This together with $\Delta_\mu e_\epsilon(y_*) \leq B_R + \dfrac{2\alpha}{(s - R)^2}(\epsilon_p + \epsilon_q)$ yields \begin{align*}\Delta_\mu V_\epsilon (y_*) &= (\Delta_\mu(g_\delta(\rho_x)) - \Delta_\mu e_\epsilon)(y_*) \\ &\geq \delta - \left(\dfrac{2\alpha}{t^2}|g_\delta'(r_x(y_*))|\epsilon_x + \dfrac{2\alpha}{(s - R)^2}(\epsilon_p + \epsilon_q)\right).\end{align*} Note that $t, s - R, r_x(y_*) > t$ are all positive and finite. Thus we can choose $\epsilon_x, \epsilon_p$ and $\epsilon_q$ sufficiently small such that $$\delta > \dfrac{2\alpha}{t^2}|g_\delta'(r_x(y_*))|\epsilon_x + \dfrac{2\alpha}{(s - R)^2}(\epsilon_p + \epsilon_q).$$ This  implies $$\Delta_\mu V_\epsilon (y_*) > 0.$$ But $V_\epsilon$ has a local maximum at $y_*$, i.e., $\nabla V_\epsilon (y_*) = 0$ and \begin{align*}\Hess V_\epsilon (y_*) \leq 0 &\implies \Delta V_\epsilon(y_*) \leq 0 \\ &\implies \Delta_\mu V_\epsilon(y_*) = \Delta V_\epsilon(y_*) - (n + 1)\left<\nabla \varphi, \underbrace{\nabla V_\epsilon(y_*)}_{=0}\right> = \Delta V_\epsilon(y_*) \leq 0.\end{align*} This contradicts $\Delta_\mu V_\epsilon(y_*) > 0$ found above. Thus $$e(x) \leq 2t + g_\delta(t), \qquad \forall \, \,  \delta > 0.$$

We now take limits in the order below. For fixed $R$ and $t > 0$, we see that $$g_\delta(t) \xrightarrow[]{\delta \downarrow 0} g_R(t) = \dfrac{B_R}{2(\lambda_R + 1)}\left(u^2 + \dfrac{2R^{\lambda_R + 1}}{\lambda_R -1}u^{1 - \lambda_R} - \dfrac{\lambda_R + 1}{\lambda_R -1}R^2\right).$$ Observe that $$\lambda_R \xrightarrow[]{R\downarrow h} A \qquad \text{and} \qquad B_R \xrightarrow[]{R \downarrow h} B_0 = \dfrac{4A}{s}.$$ Also, for a fixed $t \in (0, h)$, all expressions in $g_R(t)$ are continuous in $R$ for $t > 0$ and $A > 1$. Thus $$g_R(t) \xrightarrow[]{R \downarrow h} g_A(t) = \dfrac{B_0}{2(A + 1)}\left(t^2 + \dfrac{2h^{A + 1}}{A - 1}t^{1 - A} - \dfrac{A + 1}{A - 1}h^2\right).$$ Since $e(x) \leq 2t + g_R(t)$ for all $R \in (h, s)$, taking the limit $R \downarrow h$ gives \begin{equation}\label{eq:5.4} e(x) \leq 2t + g_A(t), \qquad 0 < t < h.\end{equation} We obtain the desired bound by finding a value $t$ that minimizes the right-hand side of \eqref{eq:5.4}. Define $F(t) = 2t + g_A(t)$. We obtain $$\begin{cases}F'(t) = 2 + \dfrac{B_0}{A + 1}\left(t - \dfrac{h^{A + 1}}{t^A}\right) \\
F''(t) = \dfrac{B_0}{A + 1}\left(1 + \dfrac{Ah^{A + 1}}{t^{A + 1}}\right) & 0 < t < h.\end{cases}$$ Since $F''(t) > 0$ on $0 < t < h$, $F$ is strictly convex on $0 < t < h$, and so there is a unique global minimum on $(0, h)$. Note that as $t \downarrow 0$, the quotient $h^{A + 1}/t^A \rightarrow + \infty$ and so $$F'(t) \xrightarrow[]{t \downarrow 0} - \infty.$$ On the other hand, $$F'(h) = 2 > 0.$$ By continuity of $F'$, there exists $t_* \in (0, h)$ such that $F'(t_*) = 0$, i.e., $t_*$ is the unique minimum. Solving for $t_*$ in $F'(t_*) = 0$, we get \begin{equation}\label{eq:5.5}
    2 + \dfrac{B_0}{A + 1}\left(t_* - \dfrac{h^{A + 1}}{t_*^A}\right) = 0 \implies 
    \dfrac{2(A+1)}{B_0}t_*^A + t_*^{A + 1} = h^{A + 1}.\end{equation} Since $t_*^{A + 1} > 0$, we get the estimate \begin{equation}\label{eq:5.6}
    t_*^A < \dfrac{B_0h^{A + 1}}{2(A + 1)} \implies t_* < \left(\dfrac{B_0h^{A + 1}}{2(A + 1)}\right)^{1/A}.
\end{equation} Multiplying equation \eqref{eq:5.5} through by $t_*^{1-A}$, we get \begin{equation}\label{eq:5.7}
    t_*^{1 -A}h^{A + 1} = \dfrac{2(A + 1)}{B_0}t_* + t_*^2
\end{equation} We substitute \eqref{eq:5.7} in $F(t_*)$ to obtain \begin{align*}
    F(t_*) &= 2t_* + \dfrac{B_0}{2(A + 1)}\left[t_*^2 + \dfrac{2}{A - 1}\left(\dfrac{2(A+1)}{B_0}t_* + t_*^2\right) - \dfrac{A + 1}{A - 1}h^2\right] \\ &= 2t_* + \left[\dfrac{B_0}{2(A + 1)} + \dfrac{B_0}{(A + 1)(A - 1)}\right]t_*^2 + \dfrac{2}{A -1}t_* - \dfrac{B_0}{2(A -1)}h^2 \\ &=  2t_* + \dfrac{B_0}{2(A -1)}t_*^2 + \dfrac{2}{A -1}t_* - \dfrac{B_0}{2(A -1)}h^2 \\ &= \underbrace{\dfrac{B_0}{2(A - 1)}}_{> 0}\underbrace{(t_*^2 - h^2)}_{< 0} + \left(\dfrac{2A}{A -1}\right)t_*.
\end{align*} Since $0 < t_* < h \implies t_*^2 - h^2 < 0$ and $\dfrac{B_0}{2(A -1)} > 0$, the first term of $F(t_*)$ is negative. Hence $$F(t_*) \leq \left(\dfrac{2A}{A - 1}\right)t_*.$$ Using the preceding estimate on $F(t_*)$ along with the estimate \eqref{eq:5.6}, we obtain \begin{align*} e(x) \leq F(t_*) &\leq \left(\dfrac{2A}{A - 1}\right)t_* \\ &<\left(\dfrac{2A}{A - 1}\right)\left(\dfrac{B_0h^{A + 1}}{2(A + 1)}\right)^{1/A} \\ &= \dfrac{2A}{A -1}\left(\dfrac{2A}{A+1}\left(\dfrac{h^{A + 1}}{s}\right)\right)^{1/A} \\ &= \left(\dfrac{2A}{A -1}\right)\left(\dfrac{2A}{A + 1}\right)^{1/A}\left(\dfrac{h^{A + 1}}{s}\right)^{1/A} \\ &= \mathcal C_A\left(\dfrac{h^{A + 1}}{s}\right)^{1/A}\end{align*} as desired. 
\end{proof}

\end{document}